\documentclass[letterpaper, 10pt, conference]{ieeeconf}
\IEEEoverridecommandlockouts 

\usepackage{amsmath,amssymb}
\usepackage{epsfig,subfigure,url}

\newcommand{\norm}[1]{\ensuremath{\left\| #1 \right\|}}
\newcommand{\bracket}[1]{\ensuremath{\left[ #1 \right]}}
\newcommand{\braces}[1]{\ensuremath{\left\{ #1 \right\}}}
\newcommand{\parenth}[1]{\ensuremath{\left( #1 \right)}}
\newcommand{\refeqn}[1]{(\ref{eqn:#1})}

\newcommand{\tr}[1]{\mbox{tr}\ensuremath{\negthickspace\bracket{#1}}}

\newcommand{\SO}{\ensuremath{\mathrm{SO(3)}}}
\newcommand{\so}{\ensuremath{\mathfrak{so}(3)}}

\renewcommand{\Re}{\ensuremath{\mathbb{R}}}

\title{\LARGE \bf
Time Optimal Attitude Control for a Rigid Body}

\author{Taeyoung Lee\authorrefmark{1}\authorrefmark{2}, Melvin Leok\authorrefmark{1}, and N. Harris McClamroch\authorrefmark{2}%
\thanks{Taeyoung Lee, and N. Harris McClamroch, Aerospace Engineering, University of Michigan, Ann Arbor, MI 48109 {\tt \{tylee,nhm\}@umich.edu}}%
\thanks{Melvin Leok, Mathematics, Purdue University, West Lafayette, IN 47907 {\tt mleok@math.purdue.edu}}%
\thanks{\textsuperscript{\footnotesize\ensuremath{*}}This research has been supported in part by NSF under grants DMS-0504747 and DMS-0726263.}
\thanks{\textsuperscript{\footnotesize\ensuremath{\dagger}}This research has been supported in part by NSF under grant CMS-0555797.}
}

\begin{document}
\allowdisplaybreaks
\maketitle \thispagestyle{empty} \pagestyle{empty}

\begin{abstract}
A time optimal attitude control problem is studied for the dynamics of a rigid body. The objective is to minimize the time to rotate the rigid body to a desired attitude and angular velocity while subject to constraints on the control input. Necessary conditions for optimality are developed directly on the special orthogonal group using rotation matrices. They completely avoid singularities associated with local parameterizations such as Euler angles, and they are expressed as compact vector equations. In addition, a discrete control method based on a geometric numerical integrator, referred to as a Lie group variational integrator, is proposed to compute the optimal control input. The computational approach is geometrically exact and numerically efficient. The proposed method is demonstrated by a large-angle maneuver for an elliptic cylinder rigid body.
\end{abstract}

\section{Introduction}
The time optimal control of spacecraft has received consistent interest as rapid attitude maneuvers are critical to various space missions such as military observation and satellite communication. The objective is to reorient the attitude of the spacecraft in a minimal maneuver time with constrained control moments. To accomplish many space missions, large-angle attitude maneuvering capabilities are required.

Time optimal attitude maneuvers have been extensively studied in the literature~\cite{ScrTho.JGCD94}. The time optimal solution is found for a single degree of freedom system, where the attitude maneuver is constrained to an eigen-axis rotation, in~\cite{Ett.AIAA89}. It is known that the eigen-axis rotation is not generally time optimal~\cite{BilWie.JGCD93,SeyKum.JGCD93}. The attitude dynamics is often simplified in an optimality analysis, e.g., by assuming an inertially symmetric rigid body model~\cite{BilWie.JGCD93,SeyKum.JGCD93,ModBha.CDC06}, linearization~\cite{BeyVad.JGCD93} and constant magnitude angular velocity~\cite{ModBha.CDC06}.

The attitude is defined as the orientation of a body-fixed frame with respect to a reference frame, and it is represented by a rotation matrix that lies on the special orthogonal group, \SO. However, most existing optimal control scheme for the dynamics of a rigid body uses coordinate representations such as Euler angles and quaternions. The minimal attitude representations like Euler angles and Rodrigues parameters have singularities, so they are not desirable for large-angle maneuvers. The non-minimal attitude representations like quaternions have associated problems. Besides the unit norm constraint, the quaternion representation double covers \SO. So, it has an inevitable ambiguity in expressing the attitude.

The objective of this paper is to solve the time optimal attitude control problem directly on $\SO$ using rotation matrices without need of any attitude parameterization. Using a specific property of the special orthogonal group, namely that the Lie algebra $\so$ is isomorphic to $\Re^3$, necessary conditions for optimality are developed and represented as vector equations on $\Re^3$. They avoid singularities associated with Euler angles completely, and the resulting expressions for the optimality necessary conditions are more compact than expressions obtained by using quaternions. Consequently, the attitude dynamics need not be simplified to make the optimal control problem tractable.

The remaining part of this paper is focused on developing a computational approach to solve this optimal control problem. The dynamics of a rigid body has certain geometric features; in addition to the configuration space being a Lie group, the dynamics are characterized by symplectic, momentum and energy preserving properties.  The most common numerical integration methods, including the widely used (non-symplectic) explicit Runge--Kutta schemes, preserve neither the Lie group structure nor these geometric properties.

Lie group variational integrators are geometric numerical integrators that preserve these geometric features of the rigid body dynamics~\cite{CMA07}. Based on this structure-preserving numerical integrator, computational approaches have been proposed to solve various optimal control problems for the dynamics of rigid bodies~\cite{CDC06.opt,ACC07.opt,CDC07.opt}. In this paper, the time optimal attitude control problem is discretized at the level of the initial problem formulation, and discrete necessary conditions for optimality are developed using the Lie group variational integrator. This provides geometrically exact but computationally efficient tools.

In summary, the optimization scheme for time optimal attitude maneuvers that we present in this paper has the following important features: (i) necessary conditions for optimality are developed directly on $\SO$, and (ii) a computational approach is adopted by using a Lie group variational integrator for overall numerical accuracy and efficiency.

This paper is organized as follows. The time optimal attitude control problem is formulated, and continuous-time necessary conditions for optimality are developed in Section II, and in a parallel fashion, a discrete-time optimal control method is presented in Section III, followed by numerical examples in Section IV.

\section{Time Optimal Attitude Control}

\subsection{Equations of Motion}
We consider the attitude dynamics of a rigid body. The configuration space is the special orthogonal group $\SO$,
\begin{align*}
    \SO=\braces{R\in\Re^{3\times 3}\,\big|\, R^T R=I_{3\times 3},\quad \det{R}=1},
\end{align*}
where the rotation matrix $R\in\SO$ represents the linear transformation from the body-fixed frame to the inertial frame.

The continuous equations of motion for the attitude dynamics of a rigid body are given by
\begin{gather}
    J\dot\Omega+\Omega\times J\Omega = u,\label{eqn:Omegadot}\\
    \dot R= R\hat\Omega,\label{eqn:Rdot}
\end{gather}
where the matrix $J\in\Re^{3\times 3}$ is the moment of inertia matrix, the vector $\Omega\in\Re^3$ is the angular velocity expressed in the body-fixed frame, and the external control moment is denoted by $u\in\Re^3$. The \textit{hat map} $\hat\cdot:\Re^3\mapsto\so$ is an isomorphism from $\Re^3$ to skew-symmetric matrices $\so$, and is defined by the condition $\hat x y = x\times y$ for all $x,y\in\Re^3$. The inverse map is denoted by the \textit{vee map} $(\cdot)^\vee:\so\mapsto\Re^3$.

\subsection{Time Optimal Attitude Control Problem}
The objective of the time optimal attitude control problem is to transfer the given initial attitude and the angular velocity $(R_\circ,\Omega_\circ)$ of the rigid body to the desired values $(R_f,\Omega_f)$ within a minimal maneuver time $t_f$ with constrained control moment $\norm{u}\leq \overline u$ for a given control limit $\overline u\in\Re$.
\begin{gather*}
\text{For given: } (R_\circ,\Omega_\circ), (R_f,\Omega_f), \bar{u}\\
\min_u \braces{\mathcal{J}=\int_{0}^{t_f} 1\,dt},\\
\text{such that } R(t_f)=R_f,\, \Omega(t_f)=\Omega_f,\\
\text{subject to } \norm{u(t)} \leq \bar{u}\;\;\forall t\in[0,t_f] \text{ and }\refeqn{Omegadot}, \refeqn{Rdot}.
\end{gather*}

\subsection{Necessary Conditions for Optimality}
We solve this optimal control problem using variational principles applied on $\SO$. Expressions for variations of a rotation matrix, and transversality conditions are presented, and necessary conditions for optimality are developed.

\paragraph*{Expressions for variations}
We represent a variation of a rotation matrix using the exponential map, $\exp:\so\mapsto\SO$
\begin{align}
    R^\epsilon = R \exp \epsilon\hat\eta,
\end{align}
where $\epsilon\in(-c,c)$ for $c>0$, and $\hat\eta\in\so$ for $\eta\in\Re^3$. Since the exponential map is a local diffeomorphism, this expression is well-defined for some constant $c$ for given $\hat\eta$. The infinitesimal variation of the rotation matrix is given by
\begin{align}
    \delta R = \frac{d}{d\epsilon}\bigg|_{\epsilon=0} R\exp
    \epsilon\hat\eta = R\hat\eta.\label{eqn:delR}
\end{align}
The infinitesimal variation of $R^T\dot R$ is obtained from \refeqn{Rdot} and \refeqn{delR} as
\begin{align}
    \delta (R^T\dot R) & = \delta R^T \dot R + R^T \delta \dot R,\nonumber\\
     & = -\eta R^T \dot R + R^T (\dot R \hat\eta + R\hat{\dot\eta}),\nonumber\\
     & = \hat{\dot\eta}+\hat\Omega\hat\eta - \hat\eta\hat\Omega,\nonumber\\
     & = (\dot\eta+ \Omega\times \eta) \widehat{\;}.\label{eqn:delOmega}
\end{align}
The variational expressions given by \refeqn{delR} and \refeqn{delOmega} are the key ingredients to developing necessary conditions for optimality for an arbitrary optimal attitude maneuver.

\paragraph*{Transversality conditions}
The differentials in the terminal attitude and the terminal angular velocity are composed of the variation for a fixed time and a term due to the terminal time variation. Since the terminal boundary conditions are fixed, we have the transversality conditions as
\begin{gather}
    \delta R(t_f)+\dot R(t_f)d t_f = R(t_f)\hat\eta(t_f) + \dot R(t_f) d t_f=0,\label{eqn:transR}\\
    \delta \Omega(t_f) +\dot\Omega(t_f) d t_f =0.\label{eqn:transOmega}
\end{gather}

\paragraph*{Necessary conditions for optimality}
Define the augmented performance index as
\begin{align*}
    \mathcal{J}_a = \int_{0}^{t_f} & 1 +\lambda^\Omega \cdot(
    u-\Omega\times J\Omega-J\dot\Omega)\\& + \lambda^R\cdot
    (\hat\Omega-R^T\dot R){^{\vee}}\,dt,
\end{align*}
where $\lambda^\Omega,\lambda^R\in\Re^3$ are Lagrange multipliers.

Using \refeqn{delOmega}, the infinitesimal variation of the augmented performance index is given by
\begin{align*}
    & \delta \mathcal J_a  = \int_{0}^{t_f}
    \lambda^\Omega\cdot(\delta u -\delta\Omega\times J\Omega - \Omega\times J\delta\Omega
    -J\delta\dot\Omega)\\
    & \qquad\qquad\quad +\lambda^R\cdot(\delta\Omega-\dot\eta-\Omega\times
    \eta)\,dt\\
    & +\big\{1 +\lambda^\Omega \cdot( u-\Omega\times J\Omega-J\dot\Omega) + \lambda^R\cdot
    (\hat\Omega-R^T\dot R)^\vee \big\} \Big|_{t_f}d t_f.
\end{align*}
Using integration by parts, we obtain
\begin{align*}
    & \delta \mathcal J_a = \int_{0}^{t_f}
    \lambda^\Omega\cdot(\delta u-\delta\Omega\times J\Omega - \Omega\times
    J\delta\Omega)+\dot\lambda^\Omega\cdot J\delta\Omega\\
    & \qquad\qquad\quad +\lambda^R\cdot(\delta\Omega-\Omega\times
    \eta)+\dot\lambda^R\cdot \eta \,dt\\
    &\hspace{10mm}-\{\lambda^\Omega\cdot J\delta\Omega+\lambda^R\cdot\eta\}\Big|^{t_f}_{0}\\
    &+\big\{1 +\lambda^\Omega \cdot( u-\Omega\times J\Omega-J\dot\Omega) + \lambda^R\cdot
    (\hat\Omega-R^T\dot R)^\vee \big\} \Big|_{t_f}d t_f.
\end{align*}
Since the initial attitude and the initial angular velocity are fixed, we have $\eta(0)=0$, $\delta\Omega(0)=0$. Substituting and rearranging, the infinitesimal variation of the augmented performance index is given by
\begin{align*}
    & \delta \mathcal J_a = \int_{0}^{t_f}
    \delta\Omega\cdot\{-J\Omega\times\lambda^\Omega-J(\lambda^\Omega\times\Omega)+J\dot\lambda^\Omega+\lambda^R\}\\
    & \qquad\qquad\quad
    +\eta\cdot\braces{\Omega\times \lambda^R +\dot\lambda^R} +\delta u \cdot \lambda^\Omega \,dt\\
    & \qquad\quad +\big\{1 +\lambda^\Omega \cdot(
    u-\Omega\times J\Omega) + \lambda^R\cdot\Omega\big\}\Big|_{t_f}d t_f.
\end{align*}
We choose multiplier equations and boundary conditions such that the expressions in all braces in the above equations are identically zero. Then, we have
\begin{align*}
    \delta \mathcal{J}_a & = \int_{0}^{t_f} \delta u \cdot \lambda^\Omega\,dt.
\end{align*}

The optimal control input $u$ must satisfy
\begin{align}
    \lambda^\Omega\cdot \delta u \geq 0,\label{eqn:pon}
\end{align}
for all admissible $\delta u$ in $t\in[0,t_f]$. If $\lambda^\Omega=0$ for a finite time period, the control input is not determined by \refeqn{pon}. Such solutions are referred to as singular arcs. Later, it is shown that there is no singular arc in this optimal control problem.

In summary, the necessary conditions for optimality are given by
\begin{list}{$\quad\bullet$}{\setlength{\leftmargin}{0pt}}%
\item Multiplier equations
\begin{gather}
    J\dot\lambda^\Omega +J(\Omega\times\lambda^\Omega)-J\Omega\times\lambda^\Omega+\lambda^R=0,\label{eqn:lamOmega}\\
    \dot\lambda^R +\Omega\times\lambda^R =0,\label{eqn:lamR}
\end{gather}
\item Optimality condition
\begin{align}
u = -\bar{u}\,(\lambda^\Omega/\norm{\lambda^\Omega}),\label{eqn:u}
\end{align}
\item Boundary and transversality conditions
\begin{gather}
(R(0),\Omega(0))=(R_\circ,\Omega_\circ),\\
(R(t_f),\Omega(t_f))= (R_f,\Omega_f),\\
\braces{1 +\lambda^\Omega \cdot(u-\Omega\times J\Omega) + \lambda^R\cdot\Omega}\Big|_{t_f}=0,\label{eqn:trans}
\end{gather}
\end{list}
Assuming that the rigid body is inertially symmetric, $J=I_{3\times 3}$, the multiplier equation \refeqn{lamOmega} is reduced to $\dot\lambda^\Omega +\lambda^R=0$.

These necessary conditions for optimality are valid for attitude maneuvers of arbitrary magnitude as they are developed by using the rotation matrix representation on $\SO$. Since the variation of the rotation matrix is expressed in terms of the Lie algebra $\so$ isomorphic to $\Re^3$, the multiplier equations are written as compact vector equations on $\Re^3$. The presented necessary conditions for optimality have neither the singularities inherent to Euler angles nor the ambiguities and redundancy associated with quaternions.

\subsection{Singular arc}\label{subsec:sa}
In this subsection, we show that singular arcs do not exist along a solution of this time optimal control problem. Suppose that there exist a singular interval, i.e. $\lambda^\Omega(t)=0$ for a finite time period in $[0,t_f]$. Then, the minimum principle given by \refeqn{pon} does not lead to a  well-defined condition for the optimal control input. Instead, the control input is determined by the requirement that the time derivative of $\lambda^\Omega$ is equal to zero.

Let the $2q$-th time derivative of $\lambda^\Omega$ be the lowest order derivative in which the control input $u$ appears explicitly with a coefficient that is not identically zero on the singular interval. Then, the integer $q$ is called the order of the singular arc~\cite{Bel.BK75}. Here, due to the special linear structure of this multiplier equation, the singular arc has infinite order. If the condition $\lambda^\Omega=\dot\lambda^\Omega=0$ is satisfied at a single point along the trajectory, $\lambda^R=\dot\lambda^R=0$, and these are satisfied identically throughout the trajectory independent of the control input. In this case, it is clear that the boundary condition \refeqn{trans} cannot be satisfied. Thus, there is no singular arc in an optimal solution.

\section{Discrete-time Time Optimal Attitude Control}
In this section, we present a computational approach, referred to as discrete optimal control of discrete Lagrangian systems~\cite{SEC07}, to solve the time optimal attitude control problem numerically. In this approach, the dynamics of the rigid body is discretized using the discrete Hamilton's principle, in order to obtain a Lie group variational integrator~\cite{CMA07}. The corresponding discrete equations of motion are imposed as dynamic constraints to be satisfied by using Lagrange multipliers, and necessary conditions for optimality, expressed as discrete equations on multipliers, are obtained.

This method yields substantial computational advantages in finding an optimal control solution. The discrete dynamics are more faithful to the continuous equations of motion, and consequently more accurate solutions to the optimal control problem are obtained. It has been shown that the discrete dynamics is more reliable even for controlled system as it computes the energy dissipation rate of controlled systems more accurately~\cite{MarWes.AN01}. In particular, the discrete flow of the Lie group variational integrator remains on $\SO$.

Optimal solutions, computed using an indirect approach, are usually sensitive to small variations of the multipliers. This causes difficulties, such as numerical ill-conditioning, when solving the necessary conditions for optimality expressed as a two-point boundary value problem. Sensitivity derivatives, computed using the discrete necessary conditions, are not corrupted by numerical dissipation caused by conventional numerical integration schemes. Thus, the proposed computational approach is more numerically robust, and the necessary conditions can be solved in a computationally efficient manner.

\subsection{Lie Group Variational Integrator}
Since the dynamics of a rigid body has the structure of a Lagrangian or Hamiltonian system, they are symplectic, momentum and energy preserving. These geometric features determine the qualitative behavior of the rigid body dynamics, and they can serve as a basis for theoretical study of rigid body dynamics.

In contrast, the most common numerical integration methods, including the widely used (non-symplectic) explicit Runge--Kutta schemes, preserve neither the Lie group structure nor these geometric properties.  Additionally, if we integrate \refeqn{Rdot} using a typical Runge--Kutta scheme, the quantity $R^T R$ inevitably drifts from the identity matrix as the simulation time increases.

In~\cite{CMA07}, Lie group variational integrators are constructed by explicitly adapting Lie group methods~\cite{IserMun.AN00} to the discrete variational principle~\cite{MarWes.AN01}. They have the desirable property that they are symplectic and momentum preserving, and they exhibit good energy behavior for an exponentially long time period. They also preserve the Lie group structure without the use of local charts, reprojection, or constraints. These geometrically exact numerical integration methods yield highly efficient and accurate computational algorithms for rigid body dynamics, and avoid singularities and ambiguities.

Using the results presented in~\cite{CMA07}, a Lie group variational integrator on $\SO$ for equations \refeqn{Omegadot}, \refeqn{Rdot} is given by
\begin{gather}
h \widehat{J\Omega_k}  = F_k J_d - J_dF_k^T,\label{eqn:findf}\\
R_{k+1} = R_k F_k,\label{eqn:Rkp}\\
J\Omega_{k+1} = F_k^T J\Omega_k + hu_{k+1},\label{eqn:Omegakp}
\end{gather}
where the subscript $k$ denotes the $k$-th step for a fixed integration step size $h\in\Re$. The matrix $J_d\in\Re^{3\times 3}$ is a nonstandard moment of inertia matrix defined by $J_d = \frac{1}{2}\mathrm{tr}[J]I_{3\times 3} - J\in\Re^{3\times 3}$. The matrix $F_k\in\SO$ denotes the relative attitude between adjacent integration steps.

For given $(R_k,x_k)$ and control input, \refeqn{findf} is solved to find $F_k$. Then $(R_{k+1},\Omega_{k+1})$ are obtained by \refeqn{Rkp} and \refeqn{Omegakp}. This yields a map $(R_k,\Omega_k)\mapsto(R_{k+1},\Omega_{k+1})$, and this process is repeated. The only implicit part is \refeqn{findf}, where the actual computation of $F_k$ is done in the Lie algebra $\so$ of dimension 3.

One of the distinct features of the Lie group variational integrator is that it preserves both the symplectic property and the Lie group structure of the rigid body dynamics. As such, it exhibits substantially improved computational accuracy and efficiency compared with other geometric integrators that preserve only one of these properties such as non-symplectic Lie group methods~\cite{CMDA07}. The symplectic property is important even in the case of controlled dynamics, since the dissipation rate of the total energy is typically computed inaccurately by non-symplectic integrators~\cite{MarWes.AN01}.

\subsection{Discrete-time Time Optimal Attitude Control Problem}
The objective is to transfer the rigid body in a prescribed way within a minimal discrete maneuver time $N$ with constrained control input.
\begin{gather*}
\text{For given: } (R_\circ,\Omega_\circ), (R_f,\Omega_f), \bar{u}\\
\min_{u_{k+1}} \braces{\mathcal{J}=\sum_{k=0}^{N-1} 1},\\
\text{such that } R_N=R_f,\, \Omega_N=\Omega_f,\\
\text{subject to } \norm{u_{k+1}} \leq \bar{u}\;\; \forall k\in[0,N\!-1] \text{ and }\refeqn{findf}\!-\!\refeqn{Omegakp}.
\end{gather*}

\subsection{Discrete-Time Necessary Conditions for Optimality}

\paragraph*{Expressions for variations}
Similar to \refeqn{delR}, the variation of rotation matrices $R_k$ and $F_k$ are expressed as
\begin{align}
    \delta R_k =  R_k\hat\eta_k,\quad \delta F_k = F_k\hat\xi_k\label{eqn:delRk}
\end{align}
for $\eta_k,\xi_k\in\Re^3$. Using this and \refeqn{Rkp}, the variation of $R_k^T R_{k+1}$ is given by
\begin{align}
\delta (R_k^T R_{k+1}) & = \delta R_k^T R_{k+1} + R_k^T \delta R_{k+1},\nonumber\\
& = -\hat\eta F_k + F_k \hat\eta_{k+1},\nonumber\\
& = F_k (-F_k^T \eta_k + \eta_{k+1})\widehat{\;}, \label{eqn:delFk}
\end{align}
where the property $\widehat{F^T x} = F^T \hat x F$ for any $x\in\Re^3$ and $F\in\SO$ is used in the last step.

Now we develop an expression for a constrained variation corresponding
\refeqn{findf}. Taking a variation of \refeqn{findf}, we obtain
\begin{align*}
h \widehat{J\delta\Omega_k} & = F_k \hat\xi_k J_d +  J_d \hat\xi_kF_k^T.
\end{align*}
Using the property, $\hat x  A+A^T\hat x=(\braces{\tr{A}I_{3\times 3}-A}x)\widehat{\;}\; $ for all $x\in\Re^3$ $A\in\Re^{3\times 3}$, the above
equation can be written as
\begin{align*}
h J\delta\hat\Omega_k & = \widehat{F_k\xi_k} F_kJ_d +  J_d F_k^T \widehat{F_k \xi_k},\\
& = (\braces{\tr{F_kJ_d}I_{3\times 3}-F_kJ_d}F_k\xi_k)\widehat{\;}.
\end{align*}
Thus, the vector $\xi_k$ is expressed in terms of $\delta\Omega_k$
\begin{align}
\xi_k = \mathcal{B}_k J\delta\Omega_k, \label{eqn:xik}
\end{align}
where $\mathcal{B}_k=hF_k^T\braces{\tr{F_kJ_d}I_{3\times 3}-F_kJ_d}^{-1}\in\Re^{3\times 3}$. This shows the relationship between $\delta\Omega_k$ and $\delta F_k$.

\paragraph*{Transversality conditions}
Similar to \refeqn{transOmega}, we choose the transversality conditions for the angular velocity as
\begin{align}
    \delta \Omega_N + (\Omega_N-\Omega_{N-1})\delta N =0.\label{eqn:dtrnsOmega}
\end{align}
The variation of the terminal attitude due to the terminal time change is expressed as
\begin{align*}
    R_N & \braces{\frac{1}{2}R_{N-1}^T (R_N-R_{N-1})+\frac{1}{2} R_N^T (R_N-R_{N-1})}\delta N\\ & = \frac{1}{2} R_N \braces{F_{N-1}-F_{N-1}^T}\delta N.
\end{align*}
This expression is chosen such that it respects the skew-symmetry of a Lie algebra $\so$ element. Using this, the transversality conditions for the attitude are given by
\begin{align}
    R_N \hat \eta_N + \frac{1}{2} R_N \braces{F_{N-1}-F_{N-1}^T}\delta N =0.\label{eqn:dtrnsR}
\end{align}

\paragraph*{Necessary conditions for optimality}

Define the augmented performance index as
\begin{align*}
\mathcal{J}_a & = \sum_{k=0}^{N-1} 1
+\lambda_k^{\Omega}\cdot\braces{-J\Omega_{k+1} + F_k^T J\Omega_k + hu_{k+1}}\nonumber\\
& +\lambda_k^{R}\cdot\frac{1}{2}\parenth{(F_k-F_k^T)^{\vee}-(R_{k}^TR_{k+1}-R_{k+1}^T R_k)^{\vee}}.
\end{align*}
Here we assume that the time step size $h$ is small so that the relative attitude rotation between adjacent integration steps is less than $\frac{\pi}{2}$, i.e. $\norm{(\mathrm{logm}F_k)^{\vee}}< \frac{\pi}{2}$. Then, $F_k$ is equal to $R_k^T R_{k+1}$ if and only if their skew parts are identical, which can be easily shown using Rodrigues' formula. Equation \refeqn{findf} is considered implicitly using a constrained variation.

Using \refeqn{delFk}, the infinitesimal variation of the augmented performance index is given by
\begin{align*}
& \delta\mathcal{J}_a = \sum_{k=0}^{N-1}
\lambda_k^{\Omega}\cdot \braces{h\delta u_{k+1}-J\delta \Omega_{k+1} + \delta
F_k^T J\Omega_k+F_k^T J\delta\Omega_k }\\
& + \lambda_k^{R}\cdot \frac{1}{2} \Big\{F_k(\xi_k+F_k^T\eta_k-\eta_{k+1})\widehat{\;}\\
& \hspace{16mm} +(\xi_k+F_k^T\eta_k-\eta_{k+1})\widehat{\;}\;F_k^T \Big\}^\vee\nonumber\\
& + \{ 1 +\lambda_{N-1}^{\Omega}\cdot\braces{-J\Omega_{N} + F_{N-1}^T J\Omega_{N-1} +
hu_{N}}\}\delta N\\
& + \lambda_{N-1}^{R}\cdot\frac{1}{2} (F_{N-1}-F_{N-1}^T)^{\vee}\delta N\\
& - \lambda_{N-1}^R\cdot\frac{1}{2} (R_{N-1}^TR_{N}-R_{N}^T R_{N-1})^{\vee}\delta N.
\end{align*}
Several algebraic manipulation steps are required here; (i) using the property
$\hat x  A+A^T\hat x=(\braces{\tr{A}I_{3\times 3}-A}x)\widehat{\;}\,$ for all
$x\in\Re^3$ and $A\in\Re^{3\times 3}$, the expression in the second braces is written as a vector form, (ii) equation \refeqn{xik} is substituted to express $\xi_k$ in terms of $\delta\Omega_k$, and (iii) using the fact that $\eta_0=0$, $\delta\Omega_0=0$, the summation indices for the variables at the $k+1$-th step are rewritten, which is the discrete analog of integration by parts. Then, we obtain
\begin{align}
\delta\mathcal{J}_a & = \sum_{k=0}^{N-1} \lambda^\Omega_k\cdot h \delta u_{k+1} \nonumber\\%
& + \sum_{k=1}^{N-1} \delta\Omega_k \cdot \Big\{-J\lambda^{\Omega}_{k-1} + J(F_k - \mathcal{B}_k^T \widehat{F_k^T J\Omega_k}) \lambda_k^{\Omega}\nonumber\\
& \hspace*{12mm}+\frac{1}{2}J\mathcal{B}_k^T (\tr{F_k}I-F_k)\lambda_k^R\Big\}\nonumber\\
&+ \sum_{k=1}^{N-1} \eta_k\cdot\Big\{\frac{1}{2}(\tr{F_{k-1}}I-F_{k-1})\lambda_{k-1}^R \nonumber\\
& \hspace*{12mm}- \frac{1}{2}F_k(\tr{F_{k}}I-F_{k})\lambda_k^R\Big\}\nonumber\\
& -\lambda_{N-1}^{\Omega}\cdot J\delta\Omega_N -\lambda_{N-1}^R\cdot\frac{1}{2}(\tr{F_{N-1}}I-F_{N-1}^T)\eta_N\nonumber\\
& +\{1+\lambda^\Omega_{N-1}\cdot\braces{-J\Omega_{N} + F_{N-1}^T J\Omega_{N-1} +
hu_{N}}\}\delta N\nonumber\\
& +\lambda_{N-1}^{R}\cdot\frac{1}{2}(F_{N-1}-F_{N-1}^T)^{\vee}\delta N\nonumber\\
& -\lambda_{N-1}^R\cdot\frac{1}{2} (R_{N-1}^TR_{N}-R_{N}^T R_{N-1})^{\vee}\delta N.\label{eqn:delJa1}
\end{align}
Substituting the transversality conditions \refeqn{dtrnsOmega} and \refeqn{dtrnsR}, all of the expressions in the last four lines of the above equation are reduced to
\begin{align}
\Big\{1& +\lambda^\Omega_{N-1} \cdot\braces{-J\Omega_{N-1} + F_{N-1}^T J\Omega_{N-1} +
hu_{N}}\nonumber\\
& +\lambda_{N-1}^{R}\cdot\frac{1}{4}\parenth{(F_{N-1})^2-(F_{N-1}^T)^2}^{\vee}\big\}\delta N.\label{eqn:delJa11}
\end{align}

We choose discrete multiplier equations such that the expressions in the first two braces in \refeqn{delJa1} are identically zero, and we choose boundary condition such that the expression given by \refeqn{delJa11} is equal to zero. Then, we have
\begin{align*}
    \delta \mathcal{J}_a & = \sum_{k=0}^{N-1} \lambda^\Omega_k\cdot h \delta u_{k+1}.
\end{align*}

The optimal control input $u_{k+1}$ must satisfy
\begin{align*}
    \lambda^\Omega_k \cdot \delta u_{k+1} \geq 0,
\end{align*}
for all admissible $\delta u_{k+1}$ and $k\in\braces{0,\cdots,N-1}$. Here, we do not show that there is no singular arc in the discrete-time optimal control problem. We assume that the result presented in Section \ref{subsec:sa} for the continuous-time case also applies to the discrete-time case.
In summary, the discrete necessary conditions for optimality are given by

\begin{list}{$\quad\bullet$}{\setlength{\leftmargin}{0pt}}%
\item Multiplier equations
\begin{gather}
\begin{aligned}
-J\lambda^{\Omega}_{k-1} + J & (F_k - \mathcal{B}_k^T \widehat{F_k^T J\Omega_k}) \lambda_k^{\Omega} \\ &+\frac{1}{2}J\mathcal{B}_k^T (\tr{F_k}I-F_k)\lambda_k^R=0,
\end{aligned}\label{eqn:lamOmegak}\\
(\tr{F_{k-1}}I-F_{k-1})\lambda_{k-1}^R - F_k(\tr{F_{k}}I-F_{k})\lambda_k^R=0.\label{eqn:lamRk}
\end{gather}
\item Optimality condition
\begin{align}
u_{k+1} =
    -\bar{u}\,(\lambda_k^\Omega/\norm{\lambda_k^\Omega})\label{eqn:ukp}
\end{align}
\item Boundary and transversality conditions
\begin{gather}
(R_0,\Omega_0)=(R_\circ,\Omega_\circ),\label{eqn:R0}\\
(R_N,\Omega_N)= (R_f,\Omega_f),\\
\begin{aligned}
1&+\lambda^\Omega_{N-1}\cdot\braces{-J\Omega_{N-1} + F_{N-1}^T J\Omega_{N-1} + hu_{N}}\\
 &+\lambda_{N-1}^{R}\cdot\frac{1}{4}\parenth{(F_{N-1})^2-(F_{N-1}^T)^2}^{\vee}=0.\label{eqn:BCN}
\end{aligned}
\end{gather}
\end{list}
In the above equations, the only implicit part is \refeqn{findf}. For a given initial condition $\{(R_0,\Omega_0),(\lambda^R_0,\lambda^\Omega_0)\}$, we solve \refeqn{findf} to obtain $F_0$, and we find the control input $u_1$ by \refeqn{ukp}. Then, $(R_1,\Omega_1)$ are obtained by \refeqn{Rkp} and \refeqn{Omegakp}. Using $\Omega_1$, we solve \refeqn{findf} to obtain $F_1$. Finally, $(\lambda^R_1,\lambda^\Omega_1)$ are obtained by \refeqn{lamRk} and \refeqn{lamOmegak}. This yields a map $\{(R_0,\Omega_0),(\lambda^R_0,\lambda^\Omega_0)\}\mapsto\{(R_1,\Omega_1),(\lambda^R_1,\lambda^\Omega_1)\}$, and this process is repeated.

The discrete necessary conditions for optimality are given by a two-point boundary value problem. This is to find the optimal discrete flow, multiplier, control input, and terminal maneuver time to satisfy the equations of motion \refeqn{findf}--\refeqn{Omegakp}, multiplier equations \refeqn{lamOmegak}, \refeqn{lamRk}, optimality condition \refeqn{ukp}, and boundary conditions \refeqn{R0}--\refeqn{BCN} simultaneously.

We use a neighboring extremal computational method~\cite{Bry.BK75}. A nominal solution satisfying all of the necessary conditions except the boundary conditions is chosen. The unspecified initial multiplier is updated so as to satisfy the specified terminal boundary conditions in the limit. This is also referred to as a shooting method. The main advantage of the neighboring extremal method is that the number of iteration variables is small. In other approaches, the initial guess of control input history or multiplier variables are iterated, so the number of optimization parameters are proportional to the number of discrete time steps.

A difficulty is that the extremal solutions are sensitive to small changes in the unspecified initial multiplier values. The nonlinearities also make it hard to construct an accurate estimate of sensitivity, and it may result in numerical ill-conditioning. By adopting a geometric numerical integrator, sensitivity derivatives along the discrete necessary conditions do not have numerical dissipation introduced by conventional numerical integration schemes. Thus, they are numerically more robust, and the necessary conditions can be solved computationally efficiently.

\section{Numerical Example}

We choose an elliptic cylinder for a rigid body model with semi-major axis $0.8\,\mathrm{m}$, semi-minor axis $0.2\,\mathrm{m}$, height $0.6\,\mathrm{m}$, mass $1,\mathrm{kg}$. The moment of inertia matrix is $J=\mathrm{diag}[0.04,\,0.19,\,0.17]\,\mathrm{kgm^2}$, and the maximum control inputs is chosen as $\overline u =0.1\,\mathrm{Nm}$.

The desired attitude maneuver is a rest-to-rest large angle rotation given by
\begin{align*}
    (R_\circ,\Omega_\circ)&=(I_{3\times 3},0)\\ (R_f,\Omega_f)&=(\exp \theta v,0),
\end{align*}
where $v=\frac{1}{\sqrt{3}}[1,\,1,\,1]\in\Re^3$, and $\theta$ is varied as $120^\circ$ and $180^\circ$.

\begin{figure}
\centerline{%
    \subfigure[Attitude maneuver]{%
        \includegraphics[width=0.46\columnwidth]{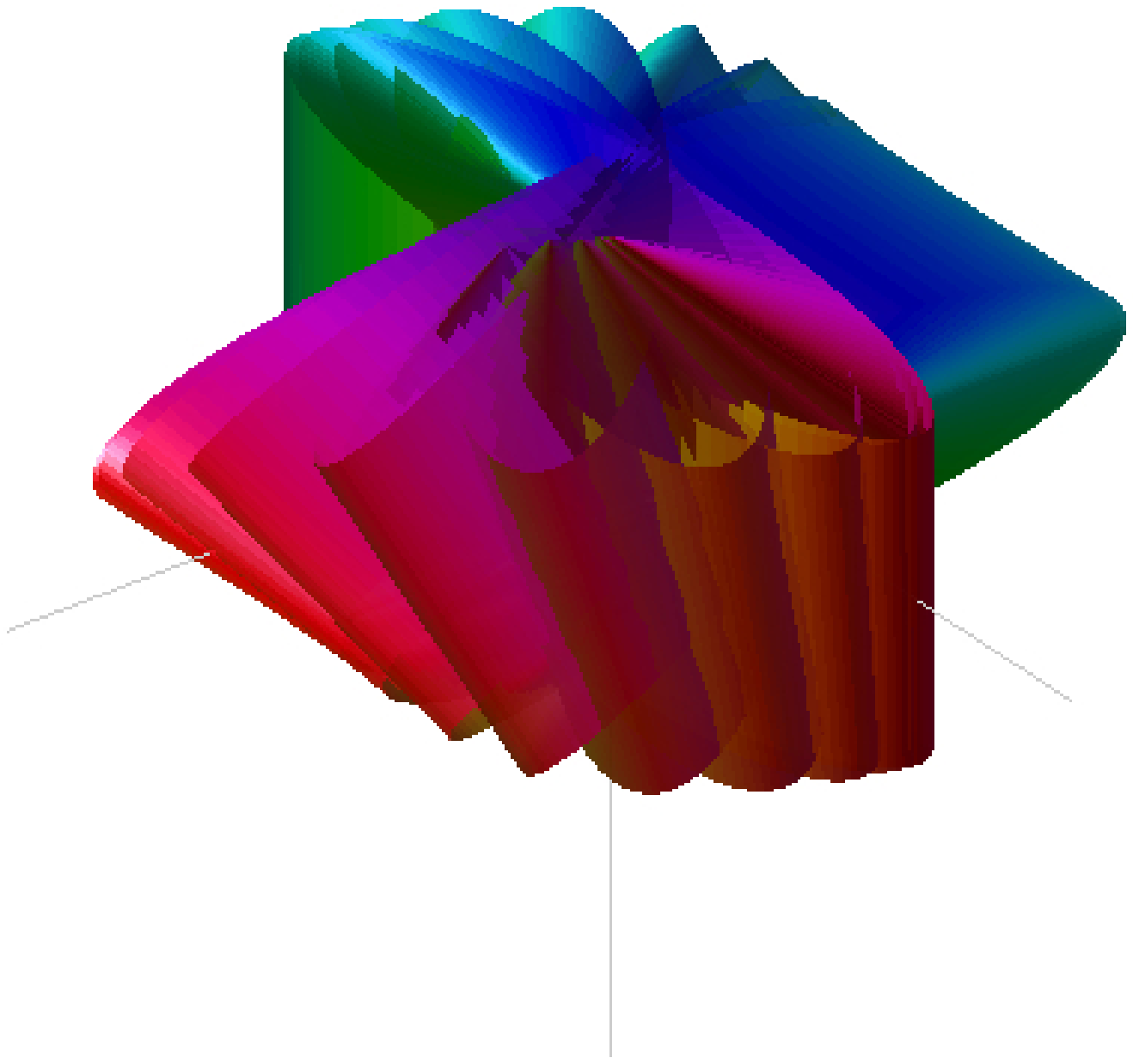}}
}
\centerline{%
    \subfigure[Angular velocity $\Omega$]{%
        \includegraphics[width=0.48\columnwidth]{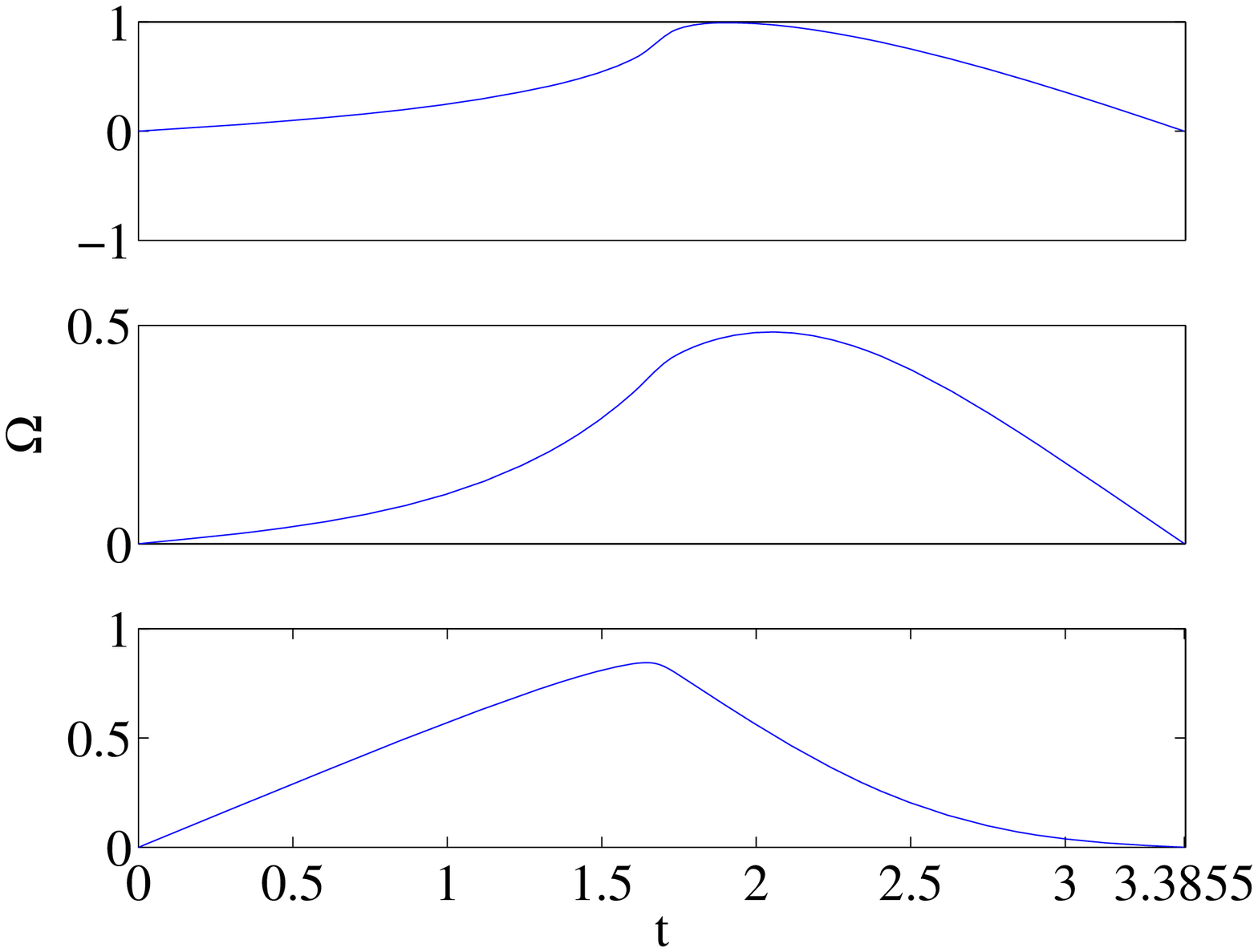}}
    \hfill
    \subfigure[Control input $u$]{%
        \includegraphics[width=0.49\columnwidth]{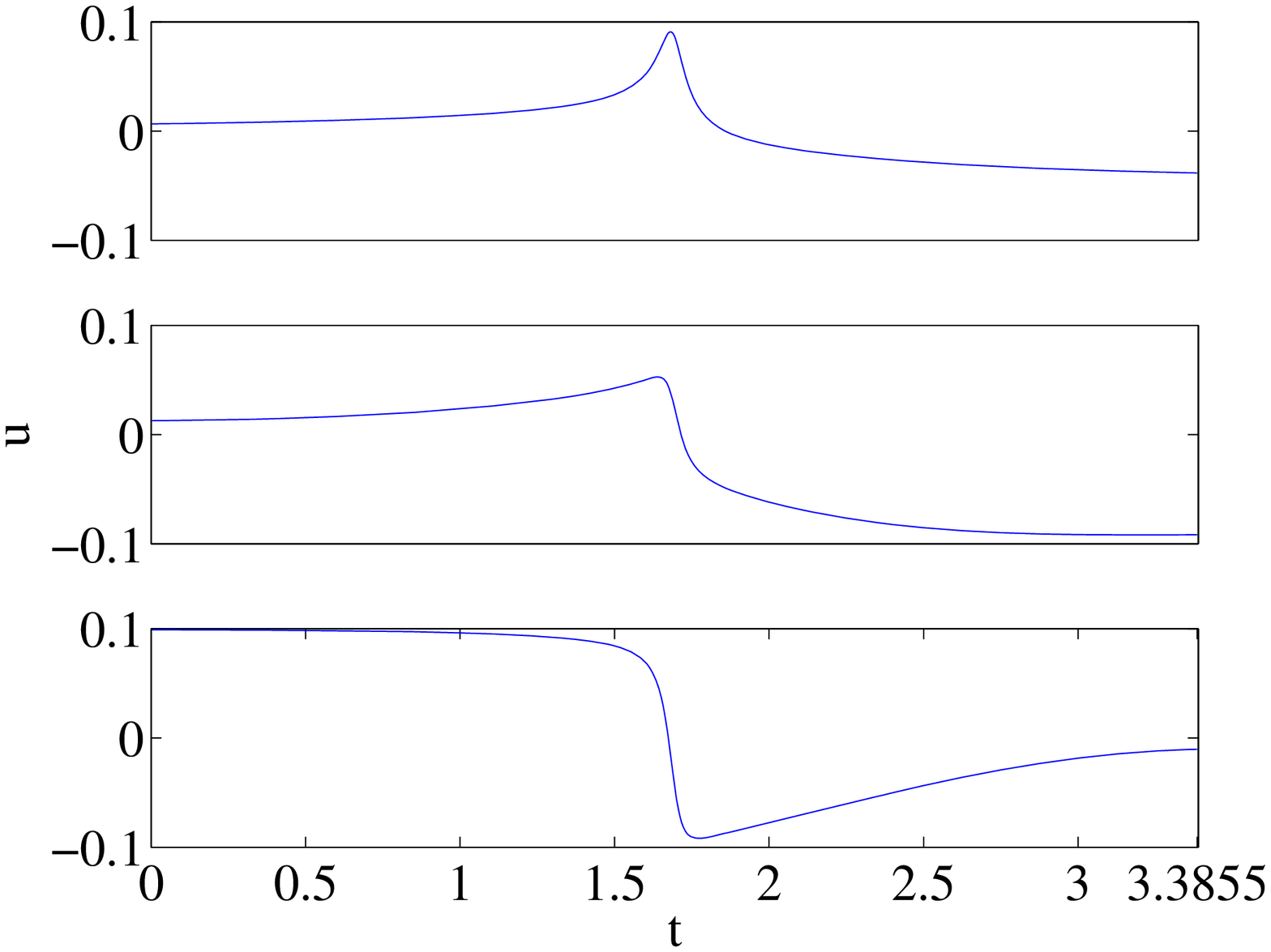}}
}
\centerline{%
    \subfigure[Lagrange multiplier $\lambda^\Omega$]{%
        \includegraphics[width=0.50\columnwidth]{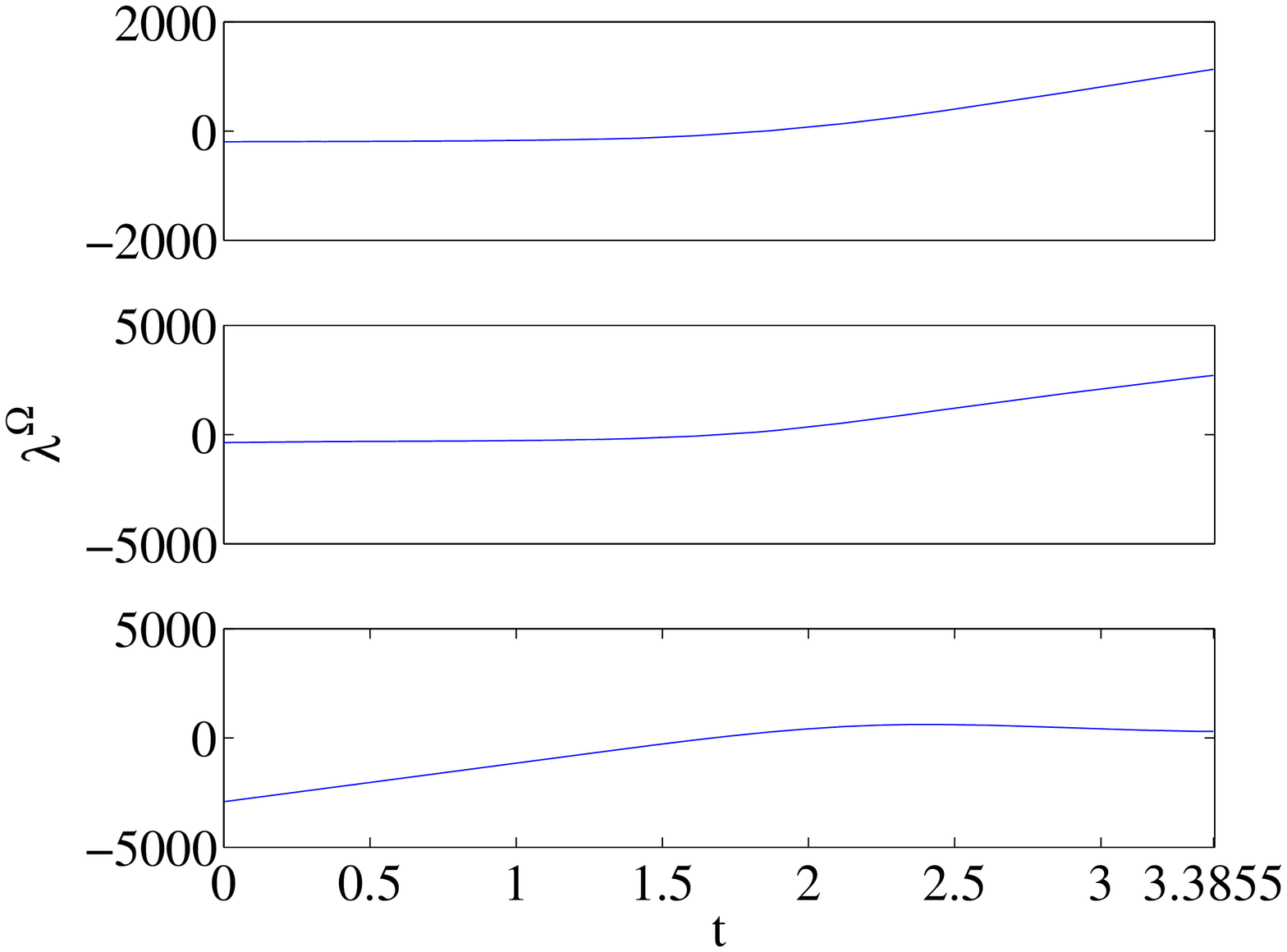}}
    \hfill
    \subfigure[Lagrange multiplier $\lambda^R$]{%
        \includegraphics[width=0.50\columnwidth]{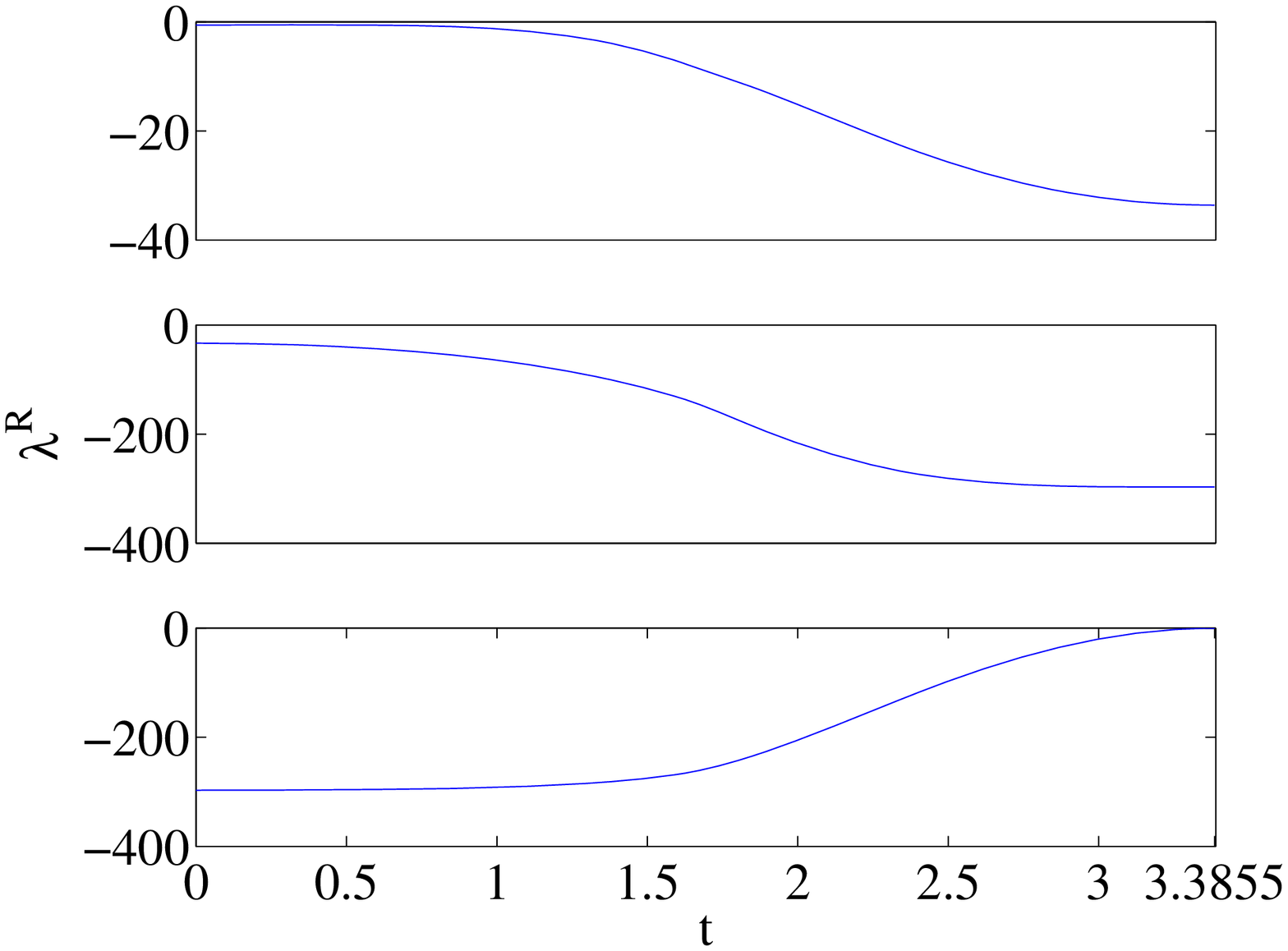}}
}
\caption{Time optimal attitude maneuver, $\theta=120^\circ$}\label{fig:d120}
\end{figure}

\begin{figure}
\centerline{%
    \subfigure[Attitude maneuver]{%
        \includegraphics[width=0.46\columnwidth]{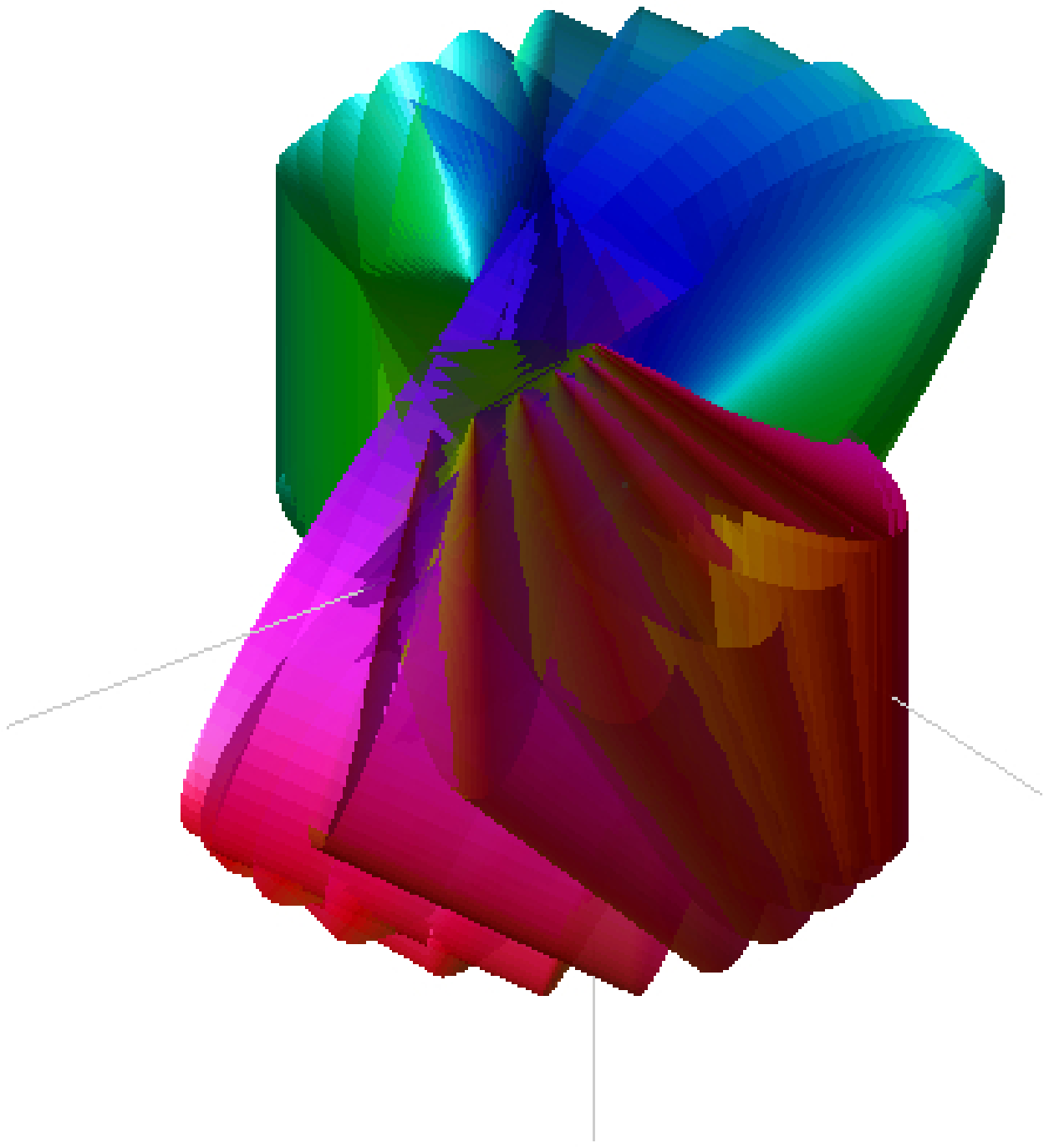}}
}
\centerline{%
    \subfigure[Angular velocity $\Omega$]{%
        \includegraphics[width=0.48\columnwidth]{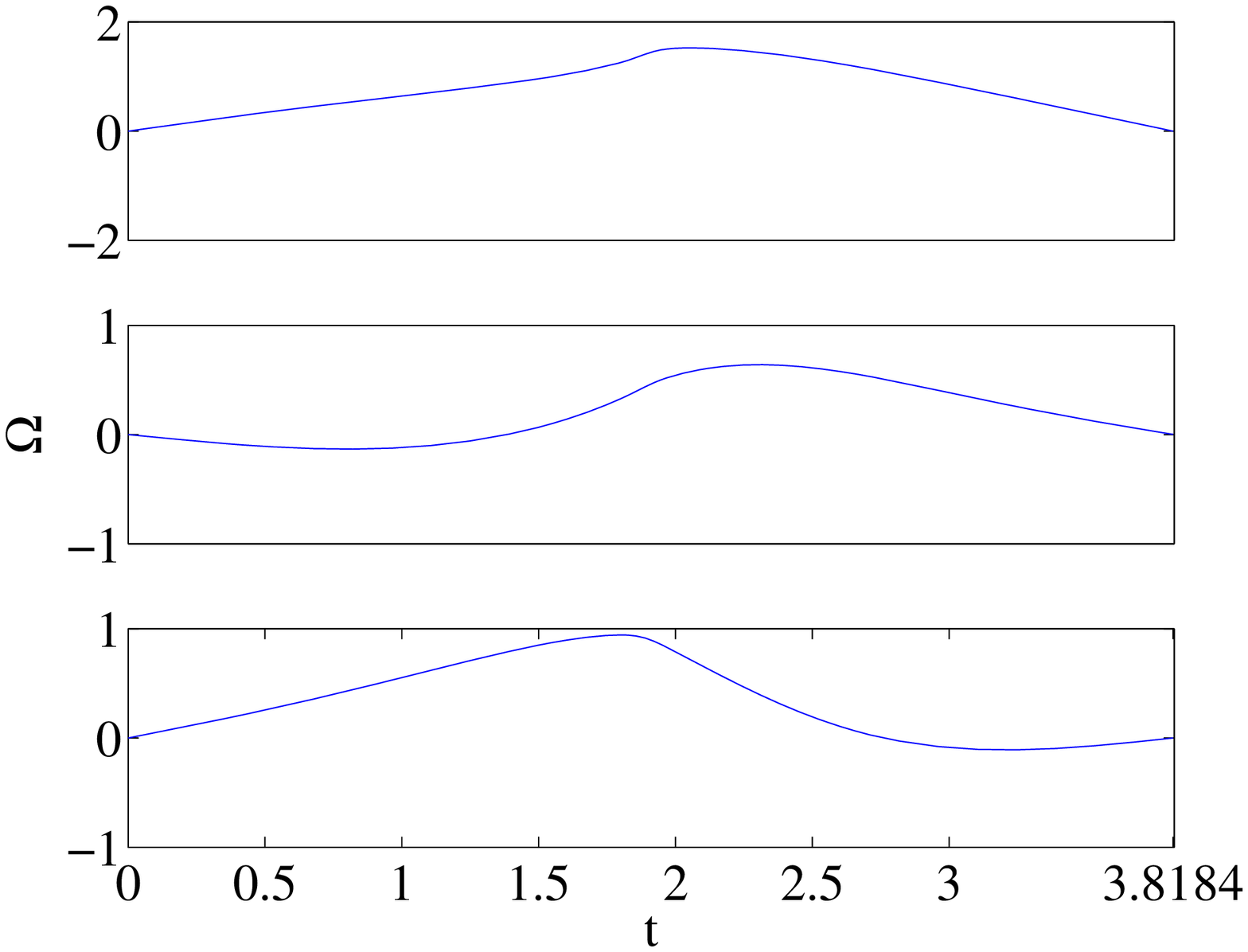}}
    \hfill
    \subfigure[Control input $u$]{%
        \includegraphics[width=0.49\columnwidth]{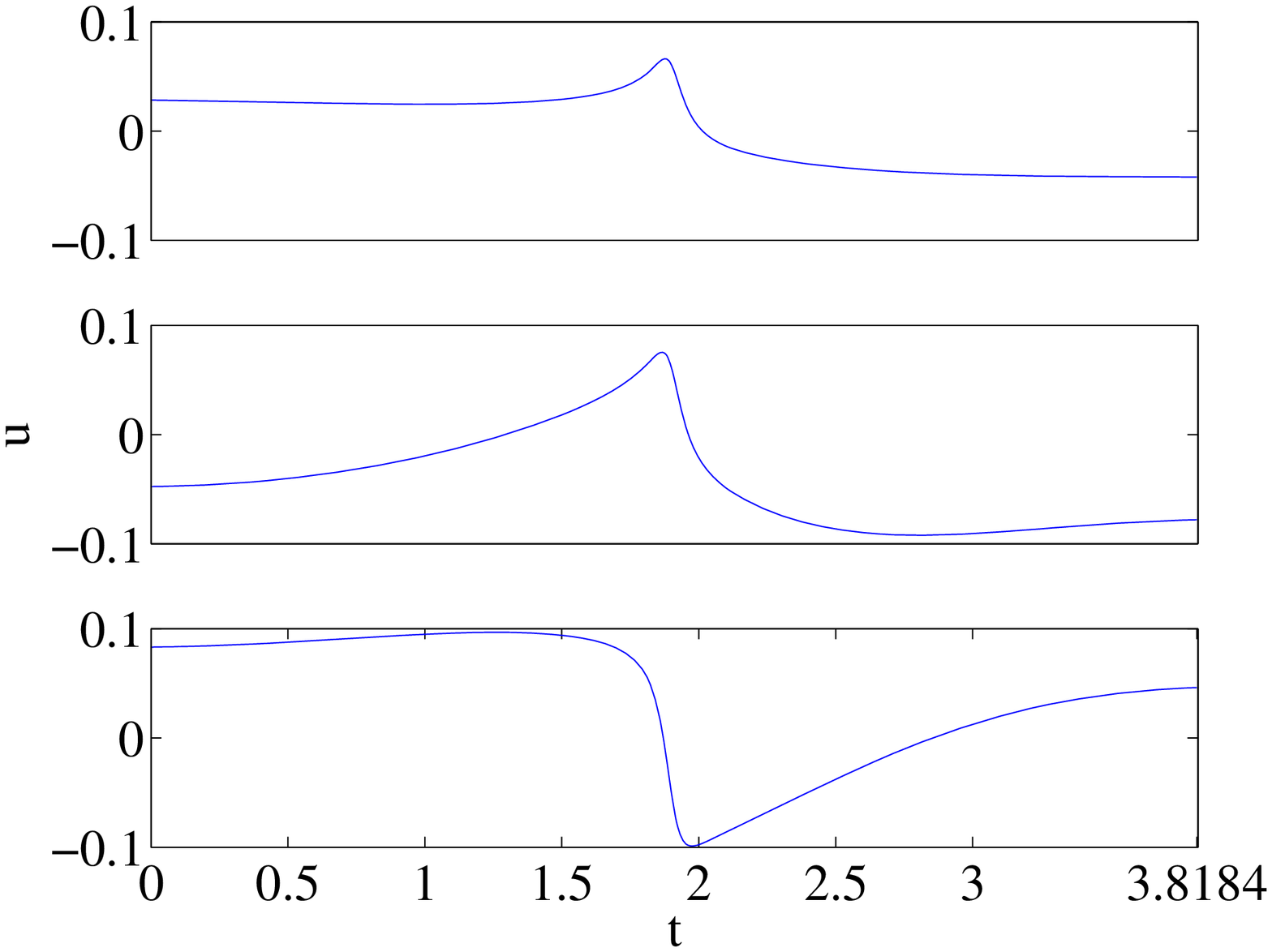}}
}
\centerline{%
    \subfigure[Lagrange multiplier $\lambda^\Omega$]{%
        \includegraphics[width=0.50\columnwidth]{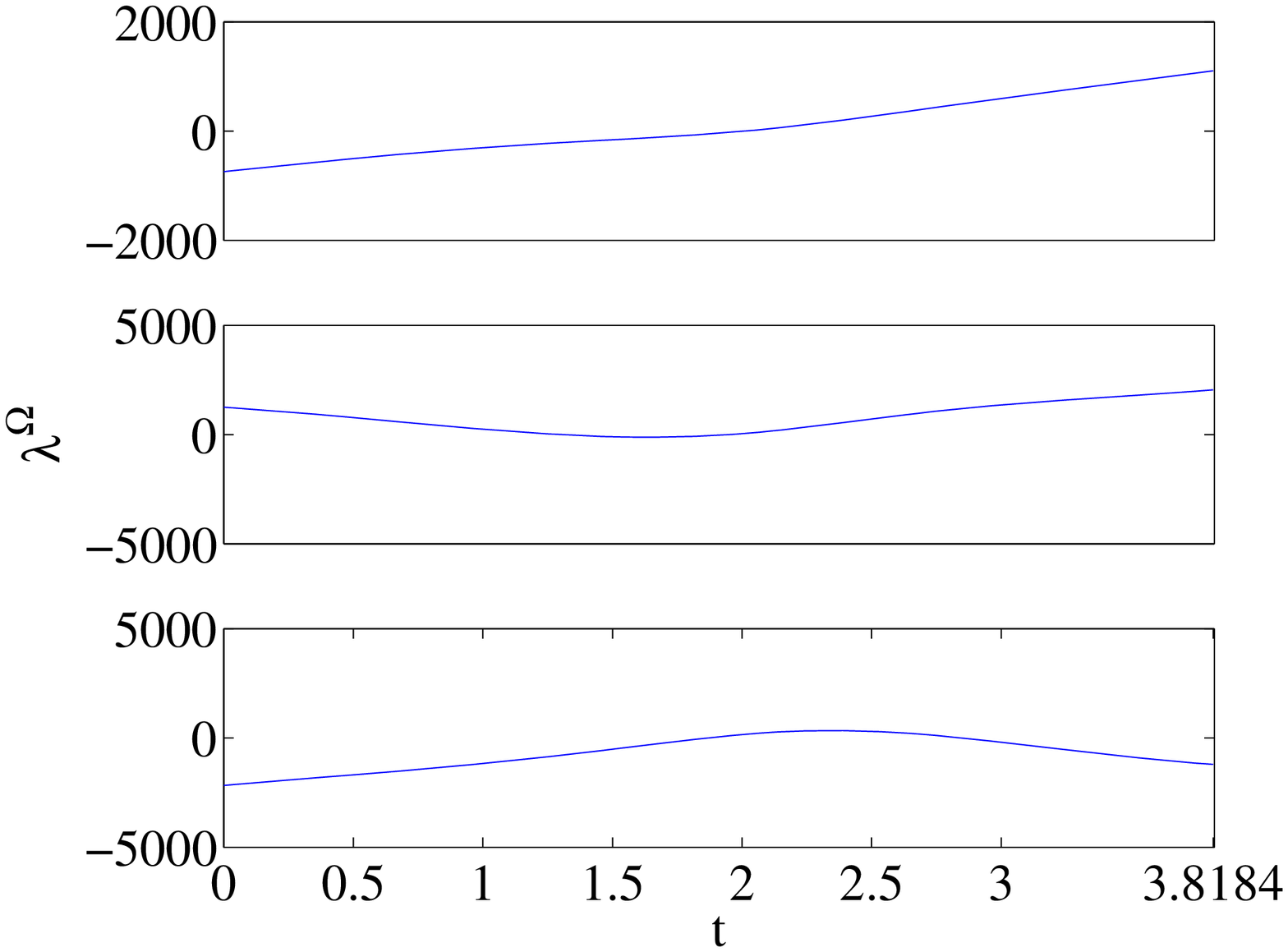}}
    \hfill
    \subfigure[Lagrange multiplier $\lambda^R$]{%
        \includegraphics[width=0.50\columnwidth]{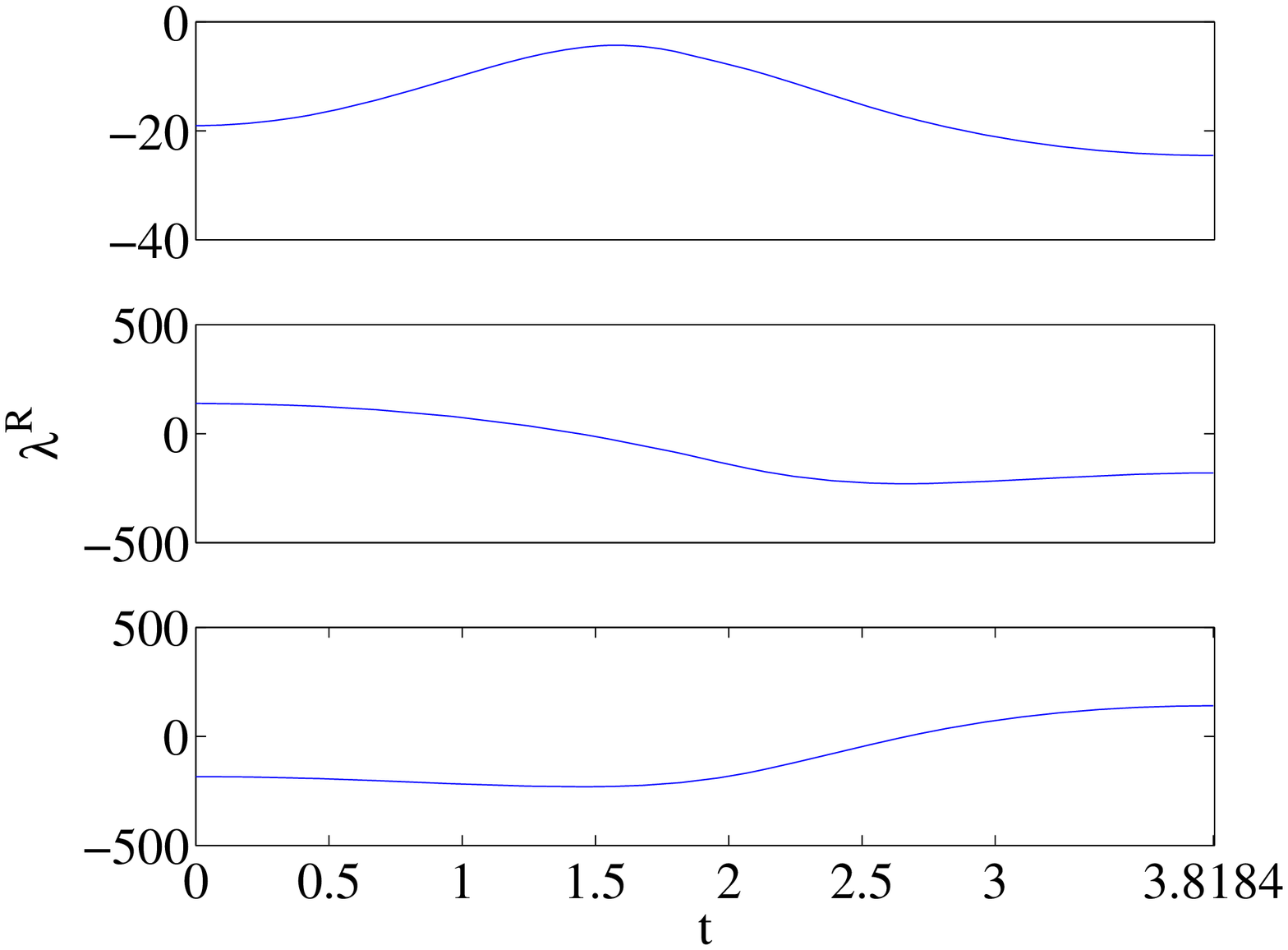}}
}
\caption{Time optimal attitude maneuver, $\theta=180^\circ$}\label{fig:d180}
\end{figure}

When deriving the discrete necessary conditions for optimality, we assume that the number of discrete steps $N$ varies. For computational purpose, it is not desirable to search the optimal value of $N$ since the terminal attitude, angular velocity and multiplier change in a discrete manner for varying integer $N$. Thus, it is not guaranteed that the boundary condition is satisfied to a desired numerical accuracy.

In the numerical computation, we fix the number of steps by an educated guess, $N=1000$ in this particular numerical example, and we vary the timestep $h$. In essence, we find the seven parameters, initial multiplier $(\lambda^R_0,\lambda^\Omega_0)$ and the time step $h$, satisfying the seven-dimensional terminal boundary conditions \refeqn{R0}--\refeqn{BCN} under the discrete equations of motion, the multiplier equation, and the optimality condition.

We solve this two-point boundary value problem, interpreted as a nonlinear equation by the shooting method, using a general nonlinear equation solver, namely the Matlab \texttt{fsolve} function. The multipliers are initialized randomly, and the timestep is initialized as $h=0.002$ seconds. The optimal solutions are found in $94$ and $211$ seconds, respectively, on Intel Pentinum M 1.73 GHz processor, and the boundary condition errors are less than $10^{-15}$.

The optimized attitude maneuver, angular velocity, multiplier, and control input history are presented in Figures \ref{fig:d120} and \ref{fig:d180}. (Simple animations which show these maneuvers of the rigid body are available at \url{http://www.umich.edu/~tylee}.) The optimized maneuver times are $3.3855$ and $3.8184$ seconds, respectively, and there is no singular arc along the optimized solutions.

\section{Conclusions}
A time optimal attitude control problem to rotate a rigid body within a minimal time with constrained control input is studied. Necessary conditions for optimality are developed on $\SO$ using rotation matrices without need of attitude parameterizations such as Euler angles and quaternions. This provides a globally applicable and compact form of necessary conditions for optimality. For overall computational accuracy and efficiency, a discrete optimal control method is proposed using a Lie group variational integrator.

In this paper, the two-norm of the control moment is constrained, and consequently, there is no singular arc in the optimal solution. The proposed necessary conditions for optimality can be directly applied, without modification, to the case where the absolute value of each component of the control moment is bounded. In this case, the expressions for optimal singular control can be developed, for example, by following the approach given in~\cite{SeyKum.JGCD93}, using the compact multiplier equations presented here.

\bibliography{opt}
\bibliographystyle{IEEEtran}

\end{document}